\documentclass[12pt, reqno]{amsart}
\usepackage{amsmath, amsthm, amscd, amsfonts, amssymb, graphicx, color}
\usepackage[bookmarksnumbered, colorlinks, plainpages]{hyperref}
\input{mathrsfs.sty}

\newtheorem{theorem}{Theorem}[section]
\newtheorem{lemma}[theorem]{Lemma}
\newtheorem{proposition}[theorem]{Proposition}
\newtheorem{corollary}[theorem]{Corollary}
\theoremstyle{definition}

\theoremstyle{remark}
\newtheorem{remark}[theorem]{Remark}
\numberwithin{equation}{section}

\begin{document}

\title[Dunkl--Williams inequality for operators]{Dunkl--Williams inequality for operators \\ associated with $p$-angular distance}

\author[F. Dadipour, M. Fujii, M.S. Moslehian]{F. Dadipour, M. Fujii and M. S. Moslehian}

\address{$^{1,3}$Department of Pure Mathematics, Center of Excellence in
Analysis on Algebraic Structures (CEAAS), Ferdowsi University of
Mashhad, P. O. Box 1159, Mashhad 91775, Iran.}
\email{dadipoor@yahoo.com}
\email{moslehian@ferdowsi.um.ac.ir and
moslehian@ams.org}
\urladdr{\url{http://profsite.um.ac.ir/~moslehian/}}
\address{$^{2}$Department of Mathematics, Osaka Kyoiku University, Kashiwara, Osaka 582-8582, Japan.}
\email{mfujii@cc.osaka-kyoiku.ac.jp}



\begin{abstract}
We present several operator versions of the Dunkl--Williams inequality
with respect to the $p$-angular distance for operators. More precisely, we show that if
$A, B \in \mathbb{B}(\mathscr{H})$ such that $|A|$ and $|B|$  are invertible, $\frac{1}{r}+\frac{1}{s}=1\,\,(r>1)$ and $p\in\mathbb{R}$, then
\begin{equation*}
|A|A|^{p-1}-B|B|^{p-1}|^{2}
\leq |A|^{p-1}(r|A-B|^{2}+s||A|^{1-p}|B|^{p}-|B||^2)|A|^{p-1}.
\end{equation*}
In the case that $0<p \leq 1$, we remove the invertibility assumption and show that if $A=U|A|$ and $B=V|B|$ are the polar decompositions of $A$ and $B$, respectively, $t>0$, then
$$|(U|A|^{p}-V|B|^{p})|A|^{1-p}|^{2}\leq (1+t)|A-B|^{2}+(1+\frac{1}{t})||B|^{p}|A|^{1-p}-|B||^2 \,.$$
We obtain several equivalent conditions, when the case of equalities hold.
\end{abstract} \maketitle


\section{Introduction}

In 1964, Dunkl and Williams \cite{D-W} showed that, for any two
nonzero vectors $x$ and $y$ in a normed space $({\mathcal
X},\|.\|)$,
\begin{eqnarray}\label{dw}
\left\| \frac{x}{\|x\|}-\frac{y}{\|y\|}\right\|
\leq\frac{4\|x-y\|}{\|x\|+\|y\|}.
\end{eqnarray}
In the same paper, the authors proved that the constant $4$ can be
replaced by $2$ if $\mathcal{X}$ is an inner product space. This
inequality has some applications in the study of geometry of Banach
spaces. Kirk and Smiley \cite{K-S} showed that inequality \eqref{dw}
with $2$ instead of $4$ characterizes inner product spaces. Thus,
the smallest number which can replace $4$ in inequality \eqref{dw}
measures ``how much'' this
space is close (or far) to be a Hilbert space, cf. \cite{J-L-M}.\\
Now the inequality \eqref{dw} is regarded as an estimation of the angular distance between given vectors $x$ and $y$.  It has many interesting refinements which have
obtained over the years, e.g. Maligranda \cite{MAL}, Merecer \cite{MER}, Dragomir \cite{DRA2}, and Pe\v{c}ari\'c and Raji\'c \cite{P-R1}.

Now we pay our attention to the following improvement of  Dunkl--Williams inequality due to Pe\v{c}ari\'c and Raji\'c:
\begin{eqnarray}\label{pr}
\left\| \frac{x}{\|x\|}-\frac{y}{\|y\|}\right\|
\leq\frac{\left(2\|x-y\|^{2}+2(\|x\|-\|y\|)^{2}\right)^\frac{1}{2}}{\max
\{ \|x\|,\|y\| \} }
\end{eqnarray}
Also they introduced an operator version of \eqref{pr} by estimating
$\left|\,A|A|^{-1}-B|B|^{-1}\,\right|$, where $A$ and $B$ are
Hilbert space operators such that $|A|$ and $|B|$ are invertible
(see Corollary \ref{c1} below).

In \cite{MAL}, Maligranda considered the $p$-angular distance ($p
\in \mathbb{R}$), as a generalization of the concept of angular
distance (when $p=0$), between nonzero elements $x$ and $y$ in a
normed space $(\mathcal{X},\|.\|)$ as
$\alpha_{p}[x,y]:= \|\|x\|^{p-1}x-\|y\|^{p-1}y \|$; see also
\cite{DRA1}.

In this paper, we introduce an operator version of the $p$-angular
distance for Hilbert space operators as a generalization of the
Pe\v{c}ari\'c--Raji\'c inequality presented in \cite{P-R3}.
Thus we will obtain the following estimation of it: If $|A|$ and $|B|$ are invertible, $\frac{1}{r}+\frac{1}{s}=1\,\,(r>1)$ and $p\in\mathbb{R}$, Then
\begin{equation*}
|A|A|^{p-1}-B|B|^{p-1}|^{2}
\leq |A|^{p-1}(r|A-B|^{2}+s||A|^{1-p}|B|^{p}-|B||^2)|A|^{p-1}.
\end{equation*}

On the other hand, Saito and Tominaga \cite{Saito-Tominaga} recently generalized Pe\v{c}ari\'c and Raji\'c inequality by deleting the invertibility condition on $|A|$ and $|B|$.  We also discuss their result.

Our basic tool is the generalized parallelogram law for operators;
$$ |A-B|^2 + \frac{1}{t} |tA+B|^2 = (1+t)|A|^2 + (1+ \frac{1}{t})|B|^2.  $$
 for any nonzero $t \in \mathbb {R}$.
We, in addition, consider several equivalent conditions when the case of equality holds in the obtained inequality. The reader is referred to \cite{FUR, MUR} for undefined notation and terminology related to Hilbert space operators.


\section{Dunkl--Williams inequality for operators}
In this section, we consider Dunkl-Williams inequality for operators as an application of generalized parallelogram law of operators (GPL):
$$ |A-B|^2 + \frac{1}{t} |tA+B|^2 = (1+t)|A|^2 + (1+ \frac{1}{t})|B|^2.  $$
for any nonzero $t \in \mathbb {R}$. This equality can be easily verified by using $|C|^2=C^*C\,\,(C \in B(H))$.

The following lemma follows from it easily.

\begin{lemma} \label{lem:lem-2.1}
Let $A,B \in \mathbb{B}(\mathscr{H})$ be operators with the polar decompositions $A=U|A|$ and $B=V|B|$. Then for each $t>0$
$$ |A-B|^2 \le (1+t)|A|^2 + (1+ \frac{1}{t})|B|^2.$$
The equality holds if and only if $tA + B = 0$.
\end{lemma}

We now state our main results, which are understood as an application of the above lemma.

\begin{theorem}\label{t0}  
Let $A,B\in \mathbb{B}(\mathscr{H})$ be operators with the polar decompositions $A=U|A|$ and $B=V|B|$ and let $t>0$
and $0 < p\leq1$ be arbitrary. Then
$$|(U|A|^{p}-V|B|^{p})|A|^{1-p}|^{2}\leq (1+t)|A-B|^{2}+(1+\frac{1}{t})||B|^{p}|A|^{1-p}-|B||^2 \,.$$
The equality holds if and only if $t(A-B) + V(|B|^{p}|A|^{1-p}-|B|) = 0$.
\end{theorem}

\begin{proof}
Replace $A$ and $B$ in the preceding lemma by $A-B$ and
$V(|B|^{p}|A|^{1-p}-|B|)$ respectively.
Then we have
\begin{eqnarray*}
 |A-V|B|^{p}|A|^{1-p})|^{2}
 &\leq&(1+t)|A-B|^{2}+(1+\frac{1}{t})|V(|B|^{p}|A|^{1-p}-|B|)|^{2} \\
 &=&(1+t)|A-B|^{2}+(1+\frac{1}{t})||B|^{p}|A|^{1-p}-|B||^{2}
\end{eqnarray*}
because $V^{*}V$ is a projection.
Hence we have the required inequality.
The equality holds if and only if $t(A-B) + V(|B|^{p}|A|^{1-p}-|B|) = 0$.
\end{proof}

Next we have an estimation of the operator $p$-angular distance.


\begin{theorem}\label{t1}
Let $A,B\in \mathbb{B}(\mathscr{H})$ such that $|A|$ and $|B|$ are invertible,
$\frac{1}{r}+\frac{1}{s}=1\,\,(r>1)$ and $p\in\mathbb{R}$. Then

\begin{equation*}
|A|A|^{p-1}-B|B|^{p-1}|^{2}
\leq|A|^{p-1}(r|A-B|^{2}+s||B|^{p}|A|^{1-p}-|B||^2|A|^{p-1}.
\end{equation*}

Moreover the equality holds if and only if
$$(r-1)(A-B)|A|^{p-1}=B(|A|^{p-1}-|B|^{p-1}).$$
\end{theorem}
\begin{proof}
The proof is similar to the above, that is, put $A_1=A-B$, $B_1= B|B|^{p-1}|A|^{1-p} - B$ and $t=r-1$ in Lemma \ref{lem:lem-2.1}.
Since $r=t+1$ and so $s=1+ \frac 1t$, we have the conclusion including the equality condition.

\end{proof}

A special case of Theorem \ref{t1}, where $p=0$ gives rise to the
main result of Pe\v{c}ari\'c and Raji\'c \cite[Theorem 2.1]{P-R3} .


\begin{corollary}\label{c1} Let $A,B\in \mathbb{B}(\mathscr{H})$ such that $|A|$ and $|B|$ are invertible and
$\frac{1}{r}+\frac{1}{s}=1\,\,(r>1)$. Then
\begin{eqnarray}\label{dad2}
|A|A|^{-1}-B|B|^{-1}|^{2}\leq
|A|^{-1}(r|A-B|^{2}+s(|A|-|B|)^{2})|A|^{-1}\,.
\end{eqnarray}
Further, the equality holds if and only if
$$(r-1)(A-B)|A|^{-1}=B(|A|^{-1}-|B|^{-1})$$
\end{corollary}
We here give some  conditions equivalent to the equality condition in
Theorem \ref{t1}.


\begin{proposition}\label{p1} Let $p\in\mathbb{R}$, $\frac{1}{r}+\frac{1}{s}=1\,\,(r>1)$ and $A,B\in \mathbb{B}(\mathscr{H})$ such that $|A|$ and $|B|$
are invertible for the case where $p<1$. Then the following
conditions are mutually equivalent:
\begin{enumerate}
\item $(r-1)(A-B)|A|^{p-1}=B(|A|^{p-1}-|B|^{p-1})\,;$
\item $(s-1)B(|A|^{p-1}-|B|^{p-1})=(A-B)|A|^{p-1}\,;$
\item $r(A-B)|A|^{p-1}+sB(|B|^{p-1}-|A|^{p-1})=0\,;$
\item $A|A|^{p-1}-B|B|^{p-1}=sB(|A|^{p-1}-|B|^{p-1})\,.$
\end{enumerate}
\end{proposition}
\begin{proof}
 The equivalence $(1)\Rightarrow(2)\Rightarrow(3)\Rightarrow(1)$ is easily checked. \\
To complete the proof, we prove
$(3)\Leftrightarrow(4)$.

Putting $t=r-1$, we have $s= \frac {t+1}t$, by which (3) and (4)  are written respectively as follows:
$$ t(A-B)|A|^{p-1}+ (t+1)B(|B|^{p-1}-|A|^{p-1})=0  $$
and
$$ tA|A|^{p-1}-B|B|^{p-1}=(t+1)B(|A|^{p-1}-|B|^{p-1}). $$
It is obvious that they are equivalent.

\end{proof}

Next we give some necessary conditions for the equality condition in Theorem \ref{t1}.

\begin{proposition}\label{t2} Let $A,B\in \mathbb{B}(\mathscr{H})$ such that $|A|$ and $|B|$ are invertible, $\frac{1}{r}+\frac{1}{s}=1\,\,(r>1)$,
$p\in\mathbb{R}$ and
\begin{eqnarray}\label{dad}
(r-1)(A-B)|A|^{p-1}=B(|A|^{p-1}-|B|^{p-1})\,. \end{eqnarray} Then
the following statements hold:
\begin{enumerate}
\item
$(r-1)|A-B|^{2}=\frac{1}{r}|A|^{1-p}|B|^{2p}|A|^{1-p}+\frac{1}{s}|A|^{2}-|B|^{2}$\,;
\item
$|B|\leq(\frac{1}{r}|A|^{1-p}|B|^{2p}|A|^{1-p}+\frac{1}{s}|A|^{2})^{\frac{1}{2}}$\,;
\item $r|A-B|=s||B|^p|A|^{1-p} - |B||$\,.
\end{enumerate}
\end{proposition}

\begin{proof}
Put $t=r-1$ and then $s= \frac {t+1}t$.

(1)
Since $t(A-B)= B(1-|B|^{p-1}|A|^{1-p})$ by the assumption, we have
$$ tA-(t+1)B= -B|B|^{p-1}|A|^{1-p}. $$
Therefore it implies that
$$ |tA-(t+1)B|^2= |A|^{1-p}|B|^{2p}|A|^{1-p}=C. $$
On the other hand, (1) is expressed  as
$$ t(t+1)|A-B|^2=C+t|A|^2-(t+1)|B|^2.  $$
So it suffices to check that
$$ |tA-(t+1)B|^2=t(t+1)|A-B|^2-t|A|^2+(t+1)|B|^2. $$

(2) It follows from (1) and the L\"owner-Heinz inequality.

(3)  Since $t(A-B)= B-B|B|^{p-1}|A|^{1-p}$ by the assumption, we have
$$ t|A-B|= |B-B|B|^{p-1}|A|^{1-p}|=||B|-|B|^p|A|^{1-p}|,  $$
which is equivalent to (3).

\end{proof}
\begin{remark}
Assume that $$(r-1)(A-B)|A|^{-1}=B(|A|^{-1}-|B|^{-1}).$$ This is the
same equation \eqref{dad} in the special case when $p=0$. From (2)
of Proposition \ref{t2} we have
$$|B|\leq (\frac{1}{r}|A|^{2}+\frac{1}{s}|A|^{2})^\frac{1}{2}=|A|$$ and so
$$\frac{r}{s}|A-B|=|A|-|B|, \text{ or }
|A|=|B|+\frac{r}{s}|A-B|,$$
 which has been shown by Pe\v{c}ari\'c and Raji\'c \cite{P-R3} .
\end{remark}

\section{Saito-Tominaga's generalization}

Very recently, Saito-Tominaga improved  Pe\v{c}ari\'c and Raji\'c inequality without the assumption of the invertibility of the absolute value of operators.

\begin{theorem}  \label{thm:Saito-Tominaga}
Let $A,B \in \mathbb{B}(\mathscr{H})$ be operators with the polar decompositions $A=U|A|$ and $B=V|B|$, and let $p,q>1$ with $\frac1p+\frac1q=1$.
Then
\begin{equation*}
|(U-V)|A||^2 \ \le \ p|A-B|^2+q(|A|-|B|)^2.
\end{equation*}
The equality
holds if and only if
\begin{equation*}
p(A-B) = qV(|B|-|A|) \quad \mbox{and} \quad V^*V = U^*U.
\end{equation*}
\end{theorem}

We here remark that it just corresponds to the case $p = 0$ in Theorem \ref{t0}.
In this section, we consider Theorem \ref{thm:Saito-Tominaga} based on the discussion in the preceding section.  For this, we rewrite it as follows:

\begin{theorem} \label{thm:Bohr45}
Let $A,B \in \mathbb{B}(\mathscr{H})$ be operators with the polar decompositions $A=U|A|$ and $B=V|B|$, and  $t>0$.
Then
\begin{equation*}
|(U-V)|A||^2 \ \le \ (t+1)|A-B|^2+ (1+ {1 \over t})(|A|-|B|)^2.
\end{equation*}The equality
holds if and only if
\begin{equation*}
t(A-B) = V(|B|-|A|) \quad \mbox{and} \quad V^*V = U^*U.
\end{equation*}
\end{theorem}

Note that Theorem \ref{thm:Saito-Tominaga} is obtained by taking $t=p-1$ in above inequality.

Now we prepare a lemma for the equality condition in above.

\begin{lemma} \label{lem:Bohr44}
Let $A,B \in \mathbb{B}(\mathscr{H})$ be operators with the polar decompositions $A=U|A|$ and $B=V|B|$ and  $t>0$.
If $t(A-B)+V(|A|-|B|)=0$ is satisfied, then
\begin{equation*}
\quad t|A-B|^2 \le |A|^2 -|B|^2,
\end{equation*}
and so $|A| \ge |B|$ and $U^*U \ge V^*V$.

In addition, if $U^*U=V^*V$, then $t|A-B|^2 = |A|^2 -|B|^2$.
\end{lemma}

\begin{proof}
Since $tA - (t+1)B = - V|A|$ by the assumption, we have
 $$|tA - (t+1)B|^2 = |A|V^*V|A|.$$
Adding $t|A|^2 - (t+1)|B|^2$ to both sides, we get
$$ t(t+1)|A-B|^2 = |A|V^*V|A| + t|A|^2 - (t+1)|B|^2 \le (t+1)(|A|^2 - |B|^2),$$
so that
$$0 \le t|A-B|^2 \le |A|^2 - |B|^2.  $$
Hence it follows that  $|A| \ge |B|$ and $U^*U \ge V^*V$.
Moreover, if $U^*U=V^*V$ is assumed, then $V^*V|A| = |A|$ and so \
$t|A-B|^2 = |A|^2 -|B|^2.$
\end{proof}

\begin{proof} {\it of Theorem \ref{thm:Bohr45}}
We replace $A$ and $B$ in
Lemma \ref{lem:lem-2.1}
by $A-B$ and $V(|A| -|B|)$ respectively.  Then we have the required inequality, and the condition for which the equality holds is that
\begin{equation*}
t(A-B) = V(|B|-|A|) \quad \mbox{and} \quad V^*V = U^*U.
\end{equation*}
The latter in above is equivalent to $|A|V^*V|A| = |A|^2$, or $V^*V|A| = |A|$, that is,
$V^*V \ge U^*U$.  By the help of the preceding
Lemma \ref{lem:Bohr44}, $|B| \le |A|$ and  $V^*V \le U^*U$, so that
$V^*V = U^*U$.
\end{proof}
\vspace{3mm}

Finally, along with the argument due to Saito and Tominaga \cite{Saito-Tominaga}, we investigate the equality condition in Theorem \ref{thm:Bohr45}.

\begin{theorem} \label{thm:Bohr48}
Let $A,B \in \mathbb{B}(\mathscr{H})$ be operators with the polar decompositions $A=U|A|$ and $B=V|B|$, and $C=W|C|$ the polar decomposition of $C=A-B$.
Assume that the equality
\begin{equation*}
|(U-V)|A||^2  = (t+1)|A-B|^2+ (1+ {1 \over t})(|A|-|B|)^2.
\end{equation*}
holds for some  $t>0$.

(1)  If $t \ge 1$, then  $A=B$.

(2)  If $0< t < 1$, then
$$ A=B(I-\frac 2{1-t}W^*W)
 \ \mbox{and} \
|A|=|B|(I+\frac {2t}{1-t}W^*W),
$$
and the converse is true.
\end{theorem}

We here prepare the following two lemmas.

\begin{lemma} \label{lem:Bohr46}
Let $A,B \in \mathbb{B}(\mathscr{H})$ be operators with the polar decompositions $A=U|A|$ and $B=V|B|$, and  $t>0$.
Suppose that $V^*V = U^*U$.
Then
$$
 t(A-B) = V(|B|-|A|)
$$
if and only if
$$ |A|=|B|+ t|A-B| \ \mbox{and} \ A-B = -V|A-B|.
$$
\end{lemma}

\begin{proof}
Since $t(A-B) = -V(|A|-|B|)$, it follows from
Lemma \ref{lem:Bohr44}
that
$$t|A-B| = ||A|-|B|| = |A| - |B|   $$
and moreover
$$ A-B= {1 \over t}V(|B|-|A|)= - {1 \over t}\ tV|A-B|=-V|A-B|. $$

Conversely, since $|A| - |B| = t|A-B|$, we have
$$t(A-B)+V(|A|-|B|)=-tV|A-B|+tV|A-B| = 0.  $$
\end{proof}
\vspace{3mm}

\begin{lemma} \label{lem:Bohr47}
Let $A,B \in \mathbb{B}(\mathscr{H})$ be operators with the polar decompositions $A=U|A|$ and $B=V|B|$, and  $t>0$.
Suppose that $V^*V = U^*U$.
If $
 t(A-B) = V(|B|-|A|)$,
then
$$ |B||A-B|+|A-B||B|=(1-t)|A-B|^2. $$
\end{lemma}
\begin{proof}
Put $C=A-B$.
The preceding lemma ensures that
$$
t|C|=|B+C|-|B|
 \ \mbox{and} \
C=-V|C|.
$$
Then it follows that
$$
|B+C|=|B|+t|C|,
$$
and that
$$
B^*C=-B^*V|C|=-(|B|V^*V)|C|=-|B||C|.
$$
Hence we have
$$ |B+C|^2 = (|B|-|C|)^2
 \ \mbox{and} \
|B+C|^2=(|B|+t|C|)^2,
$$
so that
$$
(t+1)(|B||C|+|C||B|)=(1-t^2)|C|^2,
$$
which is equivalent to the conclusion.
\end{proof}
\vspace{3mm}

Concluding this paper, we give a proof

\begin{proof} {\it of Theorem \ref{thm:Bohr48}}
The preceding lemma leads us the fact that if positive operators $S$ and $T$ satisfy $ST +TS = rS^2$ for some $r \in \mathbb{R}$, then (i) $S=0$ if $r<0$, and (ii) $S$ and $T$ commute if $r \ge 0$.  (Since $S^2T=STS-tS^3$ is selfadjoint, $S^2$ commutes with $T$ and so does $S$.)
Thus we apply it for $S=|A-B|$, $T=|B|$ and $r=1-t$.

(1) Since $r=1-t \le 0$, we first suppose that $r < 0$.  Then $S=|A-B|=0$, that is, $A=B$, as desired.  Next we suppose $r=0$.  Then $S=|C|$ commutes with $T=|B|$ and so $ST=0$.  Hence we have $|C|V^*V=0$. Moreover, since $C=-V|C|$ by Lemma \ref{lem:Bohr46}, it follows that $|C|^2= |C|V^*V|C| = 0$, i.e., $C=0$.

(2) We apply (ii).  Namely we have
$$
|B||C|=|C||B|= \frac {1-t}2 |C|^2,
$$
so that
$$
B|C|=V|B||C|= \frac {1-t}2 V|C|^2=\frac {t-1}2 C|C|=\frac {t-1}2 A|C| -\frac {t-1}2 B|C|.
$$
It implies that
$$
A|C|=\frac 2{t-1}(1+\frac {t-1}2)B|C| = \frac {t+1}{t-1}B|C|,
$$
and so
$$  AW^*W = \frac {t+1}{t-1}BW^*W.  $$
Therefore we have
$$
A=AW^*W+A(I-W^*W)=\frac {t+1}{t-1}BW^*W+B(I-W^*W)=B(I+ \frac 2{t-1}W^*W).
$$

For the second equality, it suffices to show that $W^*W$ commutes with $|B|$ because
$$|I- \frac 2{1-t}W^*W|=I+ \frac {2t}{1-t}W^*W  $$
is easily seen.  For this commutativity, we note that
$C=A-B=\frac 2{t-1}BW^*W$ by the first equality,
$C=-V|C|$ by
Lemma \ref{lem:Bohr46},
and $V^*V\ge W^*W$ by $W^*W \le \sup \{V^*V, U^*U \}$ and $V^*V=U^*U$.
So we prove that
$$
|B|W^*W=V^*BW^*W= -\frac{1-t}2 V^*C
=\frac{1-t}2 V^*V|C|=\frac{1-t}2 |C|.
$$

Incidentally the converse implication in (2) is as follows:
We first note that the second equality assures the commutativity of $|B|$ and $W^*W$.  Next it follows that
$$|A|-|B|=-\frac{2t}{1-t}|B|W^*W $$
and
$$ V|A|-B=V(|A|-|B|)=-\frac{2t}{1-t}BW^*W=-t(A-B)  $$
by the first equality.  Hence we have
$$ (U-V)|A|=A-V|A|=A+t(A-B)-B=(1+t)(A-B);   $$
$$  |(U-V)|A||^2=(1+t)^2|A-B|^2. $$
On the other hand, since
$$(|A|-|B|)^2=\left(\frac{2t}{1-t}\right)^2B^*BW^*W=t^2|A-B|^2 $$
we have
$$
(1+t)|A-B|^2 + (1+ \frac{1}{t})(|A|-|B|)^2
$$
$$
=((1+t)+((1+ \frac{1}{t})t^2)|A-B|^2 = (1+t)^2|A-B|^2.
$$
\end{proof}

\bibliographystyle{amsplain}

\end{document}